\documentclass[12pt,reqno]{amsart}
\usepackage{fullpage}
\usepackage{amsfonts}
\usepackage{amssymb}
\usepackage{times}
\usepackage{graphicx}
\newtheorem{thm}{Theorem}[section]
\newtheorem{cor}[thm]{Corollary}
\newtheorem{lem}[thm]{Lemma}

\newtheorem{Definition}[thm]{Definition}

\newtheorem{remark}[thm]{Remark}

\begin{document}

\title[COEFFICIENT ESTIMATES FOR SOME FAMILIES OF BI-BAZILEVIC FUNCTIONS OF THE MA-MINDA DEFINED BY CONVOLUTION]{\textrm{COEFFICIENT ESTIMATES FOR SOME FAMILIES OF BI-BAZILEVIC FUNCTIONS OF THE MA-MINDA DEFINED BY CONVOLUTION}}
\author{Adnan Ghazy Alamoush }

\maketitle

\begin{abstract}
\begin{flushleft}

Making use of the Hadamard product(or convolution), we find some estimates on the Taylor-Maclaurin coefficients $|a_{2}|$ and $|a_{3}|$ for functions belong to bi univalent functions of the Bazilevi$\check{c}$ type of order $\gamma$. Several (known or new) consequences of the results are also pointed out. 

\end{flushleft}
\end{abstract}
\begin{flushleft}

\end{flushleft}
\section{ Introduction}

Let $A$ denote the class of the functions $f$ of the form
\begin{equation}\label{1}
    f(z)=z+\sum^{\infty}_{n=2} a_{n}z^{n}
\end{equation}
which are analytic in the open unit disc $ U=\{z\in C:|z|<1$\} and satisfy the normalization condition $f(0)=f^{'}(0)-1=0$. Let $S$ be the subclass of $A$ consisting of functions $f$ of the form (\ref{1}) which are also univalent in $U$. Some of the important and well-investigated subclasses of the univalent function class $S$ include (for example) the class $S^{*}(\alpha)$ of starlike functions of order $\alpha$ in $U$ and the class $K(\alpha)$ of convex functions of order $\alpha$ in $U$. The well-known Koebe one-quarter theorem \cite{Duren} ensures that the image of $U$  under every univalent function $f \in A$ contains a disc of radius 1/4. Thus every univalent function $f$ has an inverse $f^{-1}$ satisfying $f^{-1}(f(z)) = z, (z \in U)$ and
$$
f^{-1}f(w)=w, (|w|<r_{0}f(z)\geq\frac{1}{4} )
$$
where
\begin{equation}\label{2}
g(w)=f^{-1}(w)=w-a_{2}w^{2}+(2a_{2}^{2}-a_{3})w^{3}-(5a_{2}^{2}-5a_{2}a_{3}+a_{4})w^{4}+...\ .
\end{equation}
We say that a function $f(z)\in A$ is  bi-univalent in $U$ if both $f(z)$ and $f^{-1}(z)$ are univalent in $U$.\\

For two functions $f$ and $g$, analytic in $U$, we say that the function $f$ is subordinate to $g$ in $U$,  written as $f(z)\prec g(x),\ \ \ (z\in U)$, provided that there exists an analytic function(that is, Schwarz function) $w(z)$ defined on $U$  with $$\ \ \  w(0)=0 \ \ and \ \  \ \ \ \ |w(z)|<1\  for\ \  all\ \ z\in U,$$
such that $f(z)=g(w(z))$ for all $z\in U$.\\

It is known that
$$f(z)\prec g(z)\ (z\in U) \ \Rightarrow\ f(0)=g(0)\ \ and\ \ f(U)\subset g(U).$$
In fact, Brannan and Taha \cite{Brannan} (see also \cite{Brannan1}) introduced certain subclasses of the bi-univalent functions similar to the familiar subclasses $S^{*}(\alpha)$ and $K(\alpha)$ of starlike and convex functions of order $ \alpha(0 \leq\alpha < 1)$, respectively, and $S_{\Sigma}^{*} (\phi)$ and $K_{\Sigma} (\phi)$ bi-starlike of Ma-Minda type and  bi-convex of Ma-Minda type,  where $\phi(z)$ is given by
\begin{equation}\label{3}
\phi(z)=1+B_{1}z+B_{2}z^{2}+B_{3}z^{3}+...\  , \ (B_{1}>0,\ z\in U).
\end{equation}
Many derivative operators can be written in terms of convolution of certain analytic functions. It is observed that this formalism brings an ease in further mathematical exploration and also helps to better understand the geometric properties of such operators. The convolution or Hadamard product of two functions $f,g\in A$ is denoted by $f * g$ and is defined as follows:
\begin{equation}\label{4}
(f*g)(z)=z+\sum^{\infty}_{n=2} a_{n}b_{n}z^{n},
\end{equation}
where $f(z)$ is given by (1) and $g(z)=z+\sum^{\infty}_{n=2} b_{n}z^{n}.$
\\

Alamoush and Darus \cite{AlAmoush} introduced differential operator $D_{\alpha,\beta,\delta,\lambda}^{k} : A \rightarrow A$ defined by
\begin{equation}\label{5}
    D_{\alpha,\beta,\delta,\lambda}^{k}f(z)=z+\sum^{\infty}_{n=2} \left[\lambda(\alpha+\beta-1)(n-1)\right]^{k}C(\delta,n)z^{n}
\end{equation}
\begin{center}
    $=z+\sum^{\infty}_{n=2}\Upsilon^{k}_{n}C(\delta,n)z^{n},$
\end{center}
for $k =0,1,2,...\ ,\  0< \alpha \leq 1,\ 0<\beta \leq1, \lambda \geq0, \delta\geq 0,\ z\in U$ , $C(\delta,n)={\begin{array}{c} \left ( \begin{array}{c} n+\delta-1\\ \delta \end{array} \right ) \end{array}},$ and $\Upsilon^{k}_{n}=\left[\lambda(\alpha+\beta-1)(n-1)\right]^{k}$.\\

Also they discussed several interesting geometrical properties exhibited by the operator $D_{\alpha,\beta,\delta,\lambda}^{k}$. Even though the parameters family of operators $D_{\alpha,\beta,\delta,\lambda}^{k}$ is a very specialized case of the widely- (and extensively-) investigated by some other authors (see \cite{Al-Oboudi}-\cite{Darus1}), tt is also easily seen that $D_{\alpha,\beta,\delta,\lambda}^{k}$ provides a generalization of the convolution between  Ruscheweyh derivative operator \cite{Ruscheweyh} and Salagean derivative operator \cite{Salagean}.\\

Recently, many authors investigated bounds for various subclasses bi-univalent function class $\Sigma$ (see\ \cite{Alamoush1}-\cite{Murugusundaramoorthy}) and obtained non-sharp coefficient estimates on the first two coefficients $|a_{2}|$ and $|a_{3}|$ of (1). For $n\ge 4$ is yet to be solved (\cite{Netanyahu}-\cite{Lewin}).\\


In 1955, Bazilevc  \cite{Bazilevc} introduced the following class of univalent functions in $U$ as:
\begin{equation}\label{6}
B_{1}(\mu)=\left\{f\in \mathbb{A}:\ \Re\left(\frac{z^{1-\mu}(f^{'}(z))}{[f(z)]^{1-\mu}}\right)> 0,\ \mu\geq0,\ z\in \mathbb{U}\right\}.
\end{equation}
This can be generalized as follows:
\begin{equation}\label{7}
B_{\alpha,\mu}=\left\{f\in \mathbb{A}:\ \Re\left(\frac{z^{1-\mu}(f^{'}(z))}{[f(z)]^{1-\mu}}\right)> \alpha,\ \  0\leq\alpha<1,\  \mu\geq0,\ z\in \mathbb{U}\right\}.
\end{equation}
Several authors have discussed various subfamilies of Bi-Bazilevi$\breve{c}$ functions of type $\gamma$. In this paper,  we find estimates on the coefficients $|a_{2}|$ and $|a_{3}|$ for functions in the new subclasses of function class $\Sigma$ involving the operator $D_{\alpha,\beta,\delta,\lambda}^{k}$. Several closely-related classes are also considered and some relevant connections to earlier known results are pointed out.

\begin{Definition}
A function $f \in \Sigma$ given by (\ref{1}) is said to be in the class $B^{k,\alpha,\beta,\delta,\lambda}_{\Sigma}(\gamma,\phi)$ if the following conditions are satisfied:
\begin{equation}\label{8}
    \frac{z^{1-\gamma}(D_{\alpha,\beta,\delta,\lambda}^{k}f(z))'}{[D_{\alpha,\beta,\delta,\lambda}^{k}f(z)]^{1-\gamma}}\prec\phi(z)
\end{equation}
and

\begin{equation}\label{9}
    \frac{z^{1-\gamma}(D_{\alpha,\beta,\delta,\lambda}^{k}g(w))'}{[D_{\alpha,\beta,\delta,\lambda}^{k}g(w)]^{1-\gamma}}\prec\phi(w),
\end{equation}
where $\gamma\geq0,\ z,w\in U$ and the function $g$ is given by (\ref{2}).
\end{Definition}

\begin{remark}
We set\ $\phi(z)=\frac{1+Az}{1-Az},\ -1\leq A<B\leq1$ then the class $B^{k,\alpha,\beta,\delta,\lambda}_{\Sigma}(\gamma,\phi)\equiv B^{k,\alpha,\beta,\delta,\lambda}_{\Sigma}(\gamma,A,B)$  which is defined as $f\in \Sigma$:
\begin{equation}\label{10}
    \frac{z^{1-\gamma}(D_{\alpha,\beta,\delta,\lambda}^{k}f(z))'}{[D_{\alpha,\beta,\delta,\lambda}^{k}f(z)]^{1-\gamma}}\prec\frac{1+Az}{1+Bz}
\end{equation}
and

\begin{equation}\label{11}
    \frac{z^{1-\gamma}(D_{\alpha,\beta,\delta,\lambda}^{k}g(w))'}{[D_{\alpha,\beta,\delta,\lambda}^{k}g(w)]^{1-\gamma}}\prec\frac{1+Aw}{1+Bw},
\end{equation}
where $z,w\in U$ and the function $g$ is given by (\ref{2}).
\end{remark}

\begin{remark}
We set\  $\phi(z)=\frac{1+(1-2\zeta)}{1-z},\ 0\leq\zeta<1$ then the class $B^{k,\alpha,\beta,\delta,\lambda}_{\Sigma}(\gamma,\phi)\equiv B^{k,\alpha,\beta,\delta,\lambda}_{\Sigma}(\gamma,\zeta)$  which is defined as $f\in \Sigma$:
\begin{equation}\label{12}
    \frac{z^{1-\gamma}(D_{\alpha,\beta,\delta,\lambda}^{k}f(z))'}{[D_{\alpha,\beta,\delta,\lambda}^{k}f(z)]^{1-\gamma}}>\zeta
\end{equation}
and

\begin{equation}\label{13}
    \frac{z^{1-\gamma}(D_{\alpha,\beta,\delta,\lambda}^{k}g(w))'}{[D_{\alpha,\beta,\delta,\lambda}^{k}g(w)]^{1-\gamma}}>\zeta,
\end{equation}
where $z,w\in U$ and the function $g$ is given by (\ref{2}).

\end{remark}
On specializing the parameters $\gamma$, one can state the various new subclasses of $\Sigma$ as illustrated in the following remarks.

\begin{remark}
For $\gamma =0$ and a function $f\in \Sigma$, given by (\ref{1}) is said to be in the class $B^{k,\alpha,\beta,\delta,\lambda}_{\Sigma}(0,\phi)=B^{k,\alpha,\beta,\delta,\lambda}_{\Sigma}(\phi)$  if the following conditions are satisfied:
\begin{equation}\label{14}
    \frac{z(D_{\alpha,\beta,\delta,\lambda}^{k}f(z))'}{D_{\alpha,\beta,\delta,\lambda}^{k}f(z)}\prec\phi(z)
\end{equation}
and

\begin{equation}\label{15}
    \frac{z(D_{\alpha,\beta,\delta,\lambda}^{k}g(w))'}{D_{\alpha,\beta,\delta,\lambda}^{k}g(w)}\prec\phi(w),
\end{equation}
where $z,w\in U$ and the function $g$ is given by (\ref{2}).
\end{remark}

\begin{remark}
For $\gamma =1$ and a function $f\in \Sigma$, given by (\ref{1}) is said to be in the class $B^{k,\alpha,\beta,\delta,\lambda}_{\Sigma}(1,\phi)=H^{k,\alpha,\beta,\delta,\lambda}_{\Sigma}(\phi)$  if the following conditions are satisfied:
\begin{equation}\label{16}
    (D_{\alpha,\beta,\delta,\lambda}^{k}f(z))'\prec\phi(z)
\end{equation}
and

\begin{equation}\label{17}
    (D_{\alpha,\beta,\delta,\lambda}^{k}g(w))'\prec\phi(w),
\end{equation}
where $z,w\in U$ and the function $g$ is given by (\ref{2}).
\end{remark}
\begin{lem}\label{lem1.2}
\cite{Pommerenke}
\textit{If $p\in P$,  then $|p_{k}|\leq2,$  for each $k$, where $P$ is the family of all functions $p$ analytic in $U$ for which $Re\{p(z)\}>0$, then}
\begin{center}
$h(z)=1+p_{1}z+p_{2}z^{2}+p_{3}z^{3}+...\ , z\in U$.
\end{center}
\end{lem}

We begin by finding the estimates on the coefficients $|a_{2}|$ and $|a_{2}|$ for functions in the class $B^{k,\alpha,\beta,\delta,\lambda}_{\Sigma}(\gamma,\phi)$.

\section{Coefficient Bounds for the Function Class  $B^{k,\alpha,\beta,\delta,\lambda}_{\Sigma}(\gamma,\phi)$ }

\begin{thm}\label{thm2.1}
\textit{
Let the function $f(z)$ given by (\ref{1}) be in the class $B^{k,\alpha,\beta,\delta,\lambda}_{\Sigma}(\gamma,\phi)$. Then}
\begin{equation}\label{18}
    |a_{2}|\leq\frac{B_{1}\sqrt{2B_{1}}}{\sqrt{\left|B_{1}^{2}\left(2(\gamma+2)\Upsilon^{k}_{3}C(\delta,3)+(\gamma-1)
    (\gamma+2)\left[\Upsilon^{k}_{2}C(\delta,2)\right]^{2}\right)-2(B_{2}-B_{1})(\gamma+1)^{2}\left[\Upsilon^{k}_{2}C(\delta,2)\right]^{2}\right|}}
\end{equation}
and
\begin{equation}\label{19}
    |a_{3}|\leq\frac{B_{1}}{(\gamma+2)\Upsilon^{k}_{3}C(\delta,3)}+\left[\frac{B_{1}}{(\gamma+1)\Upsilon^{k}_{2}C(\delta,2)}\right]^{2}.
\end{equation}
\end{thm}
{\bf Proof.} Let $f\in B^{k,\alpha,\beta,\delta,\lambda}_{\Sigma}(\gamma,\phi)$   and $g = f^{^{-1}}$. Then there are analytic functions $u,v: U\rightarrow U$ with $u(0) = 0 = v(0)$, satisfying
\begin{equation}\label{20}
    \frac{z^{1-\gamma}(D_{\alpha,\beta,\delta,\lambda}^{k}f(z))'}{[D_{\alpha,\beta,\delta,\lambda}^{k}f(z)]^{1-\gamma}}=\phi(u(z))
\end{equation}
and

\begin{equation}\label{21}
    \frac{z^{1-\gamma}(D_{\alpha,\beta,\delta,\lambda}^{k}g(w))'}{[D_{\alpha,\beta,\delta,\lambda}^{k}g(w)]^{1-\gamma}}=\phi(v(w)).
\end{equation}
Define the functions $p(z)$ and $h(z)$ by
\begin{center}
    $p(z):=\frac{1+u(z)}{1-u(z)}=1+p_{1}z+p_{2}z^{2}+p_{3}z^{3}+...\ , z\in U$,
\end{center}
and
\begin{center}
    $h(z):=\frac{1+v(z)}{1-v(z)}=1+h_{1}z+h_{2}z^{2}+h_{3}z^{3}+...\ , z\in U$.
\end{center}
or, equivalently,
\begin{equation}\label{22}
u(z):=\frac{p(z)-1}{p(z)+1}=\frac{1}{2}\left[p_{1}z\left(p_{2}-\frac{p_{1}^{2}}{2}\right)z^{2}+...\right]
\end{equation}
and
\begin{equation}\label{23}
v(z):=\frac{h(z)-1}{h(z)+1}=\frac{1}{2}\left[h_{1}z\left(h_{2}-\frac{h_{1}^{2}}{2}\right)z^{2}+...\right].
\end{equation}

Then $p(z)$ and $h(z)$ are analytic in $U$ with $p(0) = 1 = h(0)$. Since  $u, v : U\rightarrow U$ with  the
functions $p(z)$ and $h(z)$ have a positive real part in $U$, and $|p_{1}|\leq2$ and $|h_{1}|\leq2.$ Using (\ref{20}) and (\ref{21}) in (\ref{18}) and (\ref{19}) respectively, we have

\begin{equation}\label{24}
    \frac{z^{1-\gamma}(D_{\alpha,\beta,\delta,\lambda}^{k}f(z))'}{[D_{\alpha,\beta,\delta,\lambda}^{k}f(z)]^{1-\gamma}}=
    \phi\left(\frac{1}{2}\left[p_{1}z\left(p_{2}-\frac{p_{1}^{2}}{2}\right)z^{2}+...\right]\right)
\end{equation}
and
    \begin{equation}\label{25}
    \frac{z^{1-\gamma}(D_{\alpha,\beta,\delta,\lambda}^{k}g(w))'}{[D_{\alpha,\beta,\delta,\lambda}^{k}g(w)]^{1-\gamma}}=
    \phi\left(\frac{1}{2}\left[h_{1}w\left(h_{2}-\frac{h_{1}^{2}}{2}\right)w^{2}+...\right]\right).
    \end{equation}

In light of (\ref{1}) - (\ref{3}), from (\ref{22}) and (\ref{23}), it is evident that
\begin{center}
    $1+(\gamma+1)\Upsilon^{k}_{2}C(\delta,2)a_{2}z+\left[(\gamma+2)\Upsilon^{k}_{3}C(\delta,3)a_{3}+\frac{(\gamma-1)(\gamma+2)}{2}\left[\Upsilon^{k}_{2}C(\delta,2)\right]^{2}a_{2}^{2}\right]z^{2}+...$
\end{center}
\begin{center}
    $=1+\frac{1}{2}B_{1}p_{1}z+\left[\frac{1}{2}B_{1}\left(p_{2}-\frac{p_{1}^{2}}{2}\right)+\frac{1}{4}B_{2}p_{1}^{2}\right]z^{2}+...$
\end{center}
and
\begin{center}
    $1-(\gamma+1)\Upsilon^{k}_{2}C(\delta,2)a_{2}w+\left[\left(2(\gamma+2)\Upsilon^{k}_{3}C(\delta,3)+\frac{(\gamma-1)(\gamma+2)}{2}\left[\Upsilon^{k}_{2}C(\delta,2)\right]^{2}\right)a_{2}^{2}-(\gamma+2)\Upsilon^{k}_{3}C(\delta,3)a_{3}\right]w^{2}$
\end{center}

\begin{center}
    $=1+\frac{1}{2}B_{1}h_{1}w+\left[\frac{1}{2}B_{1}\left(h_{2}-\frac{h_{1}^{2}}{2}\right)+\frac{1}{4}B_{2}h_{1}^{2}\right]w^{2}+...$
\end{center}
which yield the following relations:
\begin{equation}\label{26}
    (\gamma+1)\Upsilon^{k}_{2}C(\delta,2)a_{2}=\frac{1}{2}B_{1}p_{1}
\end{equation}
\begin{equation}\label{27}
(\gamma+2)\Upsilon^{k}_{3}C(\delta,3)a_{3}+\frac{(\gamma-1)(\gamma+2)}{2}\left[\Upsilon^{k}_{2}C(\delta,2)\right]^{2}a_{2}^{2}
    =\frac{1}{2}B_{1}\left(p_{2}-\frac{p_{1}^{2}}{2}\right)+\frac{1}{4}B_{2}p_{1}^{2}
\end{equation}
\begin{equation}\label{28}
    -(\gamma+1)\Upsilon^{k}_{2}C(\delta,2)a_{2}=\frac{1}{2}B_{1}h_{1}
\end{equation}
\begin{center}
    $\left(2(\gamma+2)\Upsilon^{k}_{3}C(\delta,3)+\frac{(\gamma-1)(\gamma+2)}{2}\left[\Upsilon^{k}_{2}C(\delta,2)\right]^{2}\right)a_{2}^{2}-(\gamma+2)\Upsilon^{k}_{3}C(\delta,3)a_{3}    $
\end{center}
\begin{equation}\label{29}
=\frac{1}{2}B_{1}\left(h_{2}-\frac{h_{1}^{2}}{2}\right)+\frac{1}{4}B_{2}h_{1}^{2}.
\end{equation}
From (\ref{24}) and (\ref{26}), it follows that
\begin{equation}\label{30}
    -p_{1}=h_{1}
\end{equation}
and
\begin{equation}\label{31}
    8(\gamma+1)^{2}\left[\Upsilon^{k}_{2}C(\delta,2)\right]^{2}a_{2}^{2}=B_{1}^{2}(p_{1}^{2}+h_{1}^{2}).
\end{equation}
Adding (\ref{25}) and (\ref{27}), also making use of (\ref{28}) and (\ref{29}), we obtain
\begin{equation}\label{32}
    \left(2(\gamma+2)\Upsilon^{k}_{3}C(\delta,3)+(\gamma-1)(\gamma+2)\left[\Upsilon^{k}_{2}C(\delta,2)\right]^{2}\right)a_{2}^{2}= \frac{B_{1}}{2}(p_{2}+h_{2})+\frac{1}{4}(B_{2}-B_{1})(q_{1}^{2}+h_{1}^{2}),
\end{equation}
which yields
\begin{equation}\label{33}
    a_{2}^{2}=\frac{B_{1}^{3}(p_{2}+h_{2})}{2B_{1}^{2}\left(2(\gamma+2)\Upsilon^{k}_{3}C(\delta,3)+(\gamma-1)
    (\gamma+2)\left[\Upsilon^{k}_{2}C(\delta,2)\right]^{2}\right)-4(B_{2}-B_{1})(\gamma+1)^{2}\left[\Upsilon^{k}_{2}C(\delta,2)\right]^{2}}.
\end{equation}
Applying Lemma \ref{lem1.2} for the coefficients $p_{2}$ and $h_{2}$, we immediately have
\begin{center}

    $a_{2}^{2}\leq\frac{2B_{1}^{3}}{\left|B_{1}^{2}\left(2(\gamma+2)\Upsilon^{k}_{3}C(\delta,3)+(\gamma-1)
    (\gamma+2)\left[\Upsilon^{k}_{2}C(\delta,2)\right]^{2}\right)-2(B_{2}-B_{1})(\gamma+1)^{2}\left[\Upsilon^{k}_{2}C(\delta,2)\right]^{2}\right|}.$
\end{center}

This gives the bound on $|a_{2}|$ as asserted in (\ref{16}).
\\

Next, in order to find the bound on $|a_{2}|$, by subtracting (\ref{27}) from (\ref{25}), we get
\begin{equation}\label{34}
    2(\gamma+2)\Upsilon^{k}_{3}C(\delta,3)a_{3}-2(\gamma+2)\Upsilon^{k}_{3}C(\delta,3)a_{2}^{2}=\frac{1}{2}B_{1}\left[(p_{2}-h_{2})-\frac{1}{2}(p_{1}^{2}-h_{1}^{2})\right]
    +\frac{B_{2}}{4}(p_{1}^{2}-h_{1}^{2}).
\end{equation}
Using (\ref{28}) and (\ref{29}) in (\ref{34}), we get

\begin{center}
    $a_{3}=\frac{B_{1}(p_{2}-h_{2})}{4(\gamma+2)\Upsilon^{k}_{3}C(\delta,3)}+\frac{B_{1}^{2}\left(p_{1}^{2}+h_{1}^{2}\right)}{8(\gamma+1)
    \left[\Upsilon^{k}_{2}C(\delta,2)\right]^{2}}.$
\end{center}
\begin{flushleft}
Applying Lemma \ref{lem1.2} once again to the coefficients $p_{1}$, $q_{1}$, $p_{2}$ and $h_{2}$, we readily get (\ref{17}). This completes the proof of Theorem \ref{thm2.1}.
\end{flushleft}
Putting $\gamma = 0$ in Theorem \ref{thm2.1}, we have the following corollary.

\begin{cor}\textit{
Let the function $f(z)$ given by (\ref{1}) be in the class $B^{k,\alpha,\beta,\delta,\lambda}_{\Sigma}(\phi)$. Then
}
\begin{equation*}\label{35}
    |a_{2}|\leq\frac{B_{1}\sqrt{B_{1}}}{\sqrt{\left|B_{1}^{2}\left(2\Upsilon^{k}_{3}C(\delta,3)-\left[\Upsilon^{k}_{2}C(\delta,2)\right]^{2}\right)-(B_{2}-B_{1})\left[\Upsilon^{k}_{2}C(\delta,2)\right]^{2}\right|}}
\end{equation*}
and
\begin{equation*}\label{36}
    |a_{3}|\leq\frac{B_{1}}{2\Upsilon^{k}_{3}C(\delta,3)}+\left[\frac{B_{1}}{\Upsilon^{k}_{2}C(\delta,2)}\right]^{2}.
\end{equation*}
\end{cor}

Putting $\gamma = 1$ in Theorem \ref{thm2.1}, we have the following corollary.

\begin{cor}\textit{
Let the function $f (z)$ given by (\ref{1}) be in the class $H^{k,\alpha,\beta,\delta,\lambda}_{\Sigma}(1,\phi)$. Then
}
\begin{equation*}\label{37}
    |a_{2}|\leq\frac{B_{1}\sqrt{B_{1}}}{\sqrt{\left|3B_{1}^{2}\Upsilon^{k}_{3}C(\delta,3)-4(B_{2}-B_{1})\left[\Upsilon^{k}_{2}C(\delta,2)\right]^{2}\right|}}
\end{equation*}
and
\begin{equation*}\label{38}
    |a_{3}|\leq\frac{B_{1}}{3\Upsilon^{k}_{3}C(\delta,3)}+\left[\frac{B_{1}}{2\Upsilon^{k}_{2}C(\delta,2)}\right]^{2}.
\end{equation*}
\end{cor}
Putting $k=0$, from Corollaries 2.2 and 2.3, we get the following corollaries.
\begin{cor}\textit{
Let the function $f (z)$ given by (\ref{1}) be in the class $B^{\alpha,\beta,\delta,\lambda}_{\Sigma}(\phi)$. Then
}

\begin{equation*}\label{39}
   |a_{2}|\leq\frac{B_{1}\sqrt{B_{1}}}{\sqrt{\left|B_{1}^{2}\left(2C(\delta,3)-\left[C(\delta,2)\right]^{2}\right)-(B_{2}-B_{1})\left[C(\delta,2)\right]^{2}\right|}}
\end{equation*}
and
\begin{equation*}\label{40}
    |a_{3}|\leq\frac{B_{1}}{2C(\delta,3)}+\left[\frac{B_{1}}{C(\delta,2)}\right]^{2}.
\end{equation*}
\end{cor}
\begin{cor}\textit{
Let the function $f (z)$ given by (\ref{1}) be in the class $H^{\alpha,\beta,\delta,\lambda}_{\Sigma}(1,\phi)$. Then
}
\begin{equation*}\label{41}
    |a_{2}|\leq\frac{B_{1}\sqrt{B_{1}}}{\sqrt{\left|3B_{1}^{2}C(\delta,3)-4(B_{2}-B_{1})\left[C(\delta,2)\right]^{2}\right|}}
\end{equation*}
and
\begin{equation*}\label{42}
    |a_{3}|\leq\frac{B_{1}}{3C(\delta,3)}+\left[\frac{B_{1}}{2C(\delta,2)}\right]^{2}.
\end{equation*}
\end{cor}

Putting $\delta=0$, from Corollaries 2.2 and 2.3, we get the following corollaries.
\begin{cor}\textit{
Let the function $f (z)$ given by (\ref{1}) be in the class $B^{k,\alpha,\beta,\lambda}_{\Sigma}(\phi)$. Then
}
\begin{equation*}\label{43}
   |a_{2}|\leq\frac{B_{1}\sqrt{B_{1}}}{\sqrt{\left|B_{1}^{2}\left(2\Upsilon^{k}_{3}-[\Upsilon^{k}_{2}]^{2}\right)-(B_{2}-B_{1})[\Upsilon^{k}_{2}]^{2}\right|}}
\end{equation*}
and
\begin{equation*}\label{44}
    |a_{3}|\leq\frac{B_{1}}{2\Upsilon^{k}_{3}}+\left[\frac{B_{1}}{\Upsilon^{k}_{2}}\right]^{2}.
\end{equation*}
\end{cor}

\begin{cor}\textit{
Let the function $f (z)$ given by (\ref{1}) be in the class $H^{k,\alpha,\beta,\lambda}_{\Sigma}(1,\phi)$. Then
}
\begin{equation*}\label{45}
     |a_{2}|\leq\frac{B_{1}\sqrt{B_{1}}}{\sqrt{\left|3B_{1}^{2}\Upsilon^{k}_{3}-4(B_{2}-B_{1}[\Upsilon^{k}_{2}]^{2}\right|}}
\end{equation*}
and
\begin{equation*}\label{46}
    |a_{3}|\leq\frac{B_{1}}{3\Upsilon^{k}_{3}}+\left[\frac{B_{1}}{2\Upsilon^{k}_{2}}\right]^{2}.
\end{equation*}
\end{cor}

Putting $k=\delta=0$, from Corollaries 2.2 and 2.3, we get the following corollaries.
\begin{cor}\textit{
Let the function $f (z)$ given by (\ref{1}) be in the class $B^{\alpha,\beta,\lambda}_{\Sigma}(\phi)$. Then
}
\begin{equation*}\label{47}
     |a_{2}|\leq\frac{B_{1}\sqrt{B_{1}}}{\sqrt{\left|B_{1}^{2}-(B_{2}-B_{1})\right|}}
\end{equation*}
and
\begin{equation*}\label{48}
    |a_{3}|\leq\frac{B_{1}}{2}+B^{2}_{1}.
\end{equation*}
\end{cor}
\begin{cor}\textit{
Let the function $f(z)$ given by (\ref{1}) be in the class  $H^{\alpha,\beta,\lambda}_{\Sigma}(1,\phi)$. Then
}
\begin{equation*}\label{49}
     |a_{2}|\leq\frac{B_{1}\sqrt{B_{1}}}{\sqrt{\left|3B_{1}^{2}-4(B_{2}-B_{1})\right|}}
\end{equation*}
and
\begin{equation*}\label{50}
    |a_{3}|\leq\frac{B_{1}}{3}+\left[\frac{B_{1}}{2}\right]^{2}.
\end{equation*}
\end{cor}

\begin{remark}
Putting $k=\delta=0$ in  Theorem \ref{thm2.1}, it  reduces to Theorem 2.8 of Deniz \cite{Deniz}.
\end{remark}

\section{Corollaries and Consequences}
By setting $\phi(z)=\frac{1+Az}{1+Bz},\ -1\leq B<A\leq1 $ in Theorem \ref{thm2.1}, we state the following theorem:
\begin{thm}\label{thm3.1}
\textit{
Let the function $f(z)$ given by (\ref{1}) be in the class $B^{k,\alpha,\beta,\delta,\lambda}_{\Sigma}(\gamma,A,B)$. Then}
\begin{equation*}\label{51}
     |a_{2}|\leq\frac{\sqrt{2}(A-B)}{\sqrt{\left|(A-B)\left(2(\gamma+2)\Upsilon^{k}_{3}C(\delta,3)+(\gamma-1)
    (\gamma+2)\left[\Upsilon^{k}_{2}C(\delta,2)\right]^{2}\right)-2(B+1)(\gamma+1)^{2}\left[\Upsilon^{k}_{2}C(\delta,2)\right]^{2}\right|}}
\end{equation*}
and
\begin{equation*}\label{52}
    |a_{3}|\leq\frac{A-B}{(\gamma+2)\Upsilon^{k}_{3}C(\delta,3)}+\left[\frac{(A-B)}{(\gamma+1)\Upsilon^{k}_{2}C(\delta,2)}\right]^{2}.
\end{equation*}

\end{thm}

Putting $\gamma = 0$ in Theorem \ref{thm3.1}, we have the following corollary.

\begin{cor}
\textit{
Let the function $f(z)$ given by (\ref{1}) be in the class $B^{k,\alpha,\beta,\delta,\lambda}_{\Sigma}(A,B)$. Then}
\begin{equation*}\label{53}
     |a_{2}|\leq\frac{A-B}{\sqrt{\left|(A-B)\left(2\Upsilon^{k}_{3}C(\delta,3)-\left[\Upsilon^{k}_{2}C(\delta,2)\right]^{2}\right)-(B+1)\left[\Upsilon^{k}_{2}C(\delta,2)\right]^{2}\right|}}
\end{equation*}
and
\begin{equation*}\label{54}
    |a_{3}|\leq\frac{A-B}{2\Upsilon^{k}_{3}C(\delta,3)}+\left[\frac{(A-B)}{\Upsilon^{k}_{2}C(\delta,2)}\right]^{2}.
\end{equation*}

\end{cor}

Putting $\gamma = 1$ in Theorem \ref{thm3.1}, we have the following corollary.

\begin{cor}
\textit{
Let the function $f(z)$ given by (\ref{1}) be in the class $B^{k,\alpha,\beta,\delta,\lambda}_{\Sigma}(1,A,B)$. Then}

\begin{equation*}\label{55}
     |a_{2}|\leq\frac{A-B}{\sqrt{\left|3(A-B)\Upsilon^{k}_{3}C(\delta,3)-4(B+1)\left[\Upsilon^{k}_{2}C(\delta,2)\right]^{2}\right|}}
\end{equation*}
and
\begin{equation*}\label{56}
    |a_{3}|\leq\frac{A-B}{3\Upsilon^{k}_{3}C(\delta,3)}+\left[\frac{A-B}{2\Upsilon^{k}_{2}C(\delta,2)}\right]^{2}.
\end{equation*}

\end{cor}
Putting $k = 0$, from Corollaries 3.2 and 3.3, we get the following corollaries.

\begin{cor}
\textit{
Let the function $f (z)$ given by (\ref{1}) be in the class $B^{\alpha,\beta,\delta,\lambda}_{\Sigma}(A,B)$. Then
}

\begin{equation*}\label{57}
     |a_{2}|\leq\frac{A-B}{\sqrt{\left|(A-B)\left(2C(\delta,3)-\left[C(\delta,2)\right]^{2}\right)-(B+1)\left[C(\delta,2)\right]^{2}\right|}}
\end{equation*}
and
\begin{equation*}\label{58}
    |a_{3}|\leq\frac{A-B}{2C(\delta,3)}+\left[\frac{A-B}{C(\delta,2)}\right]^{2}.
\end{equation*}
\end{cor}
\begin{cor}\textit{
Let the function $f (z)$ given by (\ref{1}) be in the class $H^{\alpha,\beta,\delta,\lambda}_{\Sigma}(1,A,B)$. Then
}
\begin{equation*}\label{59}
     |a_{2}|\leq\frac{A-B}{\sqrt{\left|3(A-B)C(\delta,3)-4(B+1)\left[C(\delta,2)\right]^{2}\right|}}
\end{equation*}
and
\begin{equation*}\label{60}
    |a_{3}|\leq\frac{A-B}{3C(\delta,3)}+\left[\frac{A-B}{2C(\delta,2)}\right]^{2}.
\end{equation*}
\end{cor}

Putting $\delta=0$, from Corollaries 3.2 and 3.3, we get the following corollaries.
\begin{cor}

\textit{
Let the function $f(z)$ given by (\ref{1}) be in the class $B^{k,\alpha,\beta,\delta,\lambda}_{\Sigma}(A,B)$. Then}
\begin{equation*}\label{61}
     |a_{2}|\leq\frac{A-B}{\sqrt{\left|(A-B)\left(2\Upsilon^{k}_{3}-[\Upsilon^{k}_{2}]^{2}\right)-(B+1)[\Upsilon^{k}_{2}]^{2}\right|}}
\end{equation*}
and
\begin{equation*}\label{62}
    |a_{3}|\leq\frac{A-B}{2\Upsilon^{k}_{3}}+\left[\frac{(A-B)}{\Upsilon^{k}_{2}}\right]^{2}.
\end{equation*}
\end{cor}
\begin{cor}\textit{
Let the function $f (z)$ given by (\ref{1}) be in the class $H^{k,\alpha,\beta,\lambda}_{\Sigma}(1,A,B)$. Then
}
\begin{equation*}\label{63}
     |a_{2}|\leq\frac{A-B}{\sqrt{\left|3(A-B)\Upsilon^{k}_{3}-4(B+1)[\Upsilon^{k}_{2}]^{2}\right|}}
\end{equation*}
and
\begin{equation*}\label{64}
    |a_{3}|\leq\frac{A-B}{3\Upsilon^{k}_{3}}+\left[\frac{A-B}{2\Upsilon^{k}_{2}}\right]^{2}.
\end{equation*}

\end{cor}

Putting $k=\delta=0$, from Corollaries 3.2 and 3.3, we get the following corollaries.
\begin{cor}\textit{
Let the function $f (z)$ given by (\ref{1}) be in the class $B^{\alpha,\beta,\lambda}_{\Sigma}(A,B)$. Then
}
\begin{equation*}\label{65}
     |a_{2}|\leq\frac{A-B}{\sqrt{\left|(A-B)-(B+1)\right|}}
\end{equation*}
and
\begin{equation*}\label{66}
    |a_{3}|\leq\frac{A-B}{2}+(A-B)^{2}.
\end{equation*}
\end{cor}
\begin{cor}\textit{
Let the function $f (z)$ given by (\ref{1}) be in the class  $H^{\alpha,\beta,\lambda}_{\Sigma}(1,A,B)$. Then
}
\begin{equation*}\label{67}
     |a_{2}|\leq\frac{A-B}{\sqrt{\left|3(A-B)^{2}-4(B+1)\right|}}
\end{equation*}
and
\begin{equation*}\label{68}
    |a_{3}|\leq\frac{A-B}{3}+\left[\frac{A-B}{2}\right]^{2}.
\end{equation*}
\end{cor}

Further, by setting $\phi(z) = 1+\frac{1+(1-2\zeta)z}{1-z},\  0 \leq\zeta < 1$ in Theorem \ref{thm2.1} we get the following result.

\begin{thm}\label{thm3.2}
\textit{
Let the function $f(z)$ given by (\ref{1}) be in the class $B^{k,\alpha,\beta,\delta,\lambda}_{\Sigma}(\gamma,\zeta)$. Then}
\begin{equation*}\label{69}
     |a_{2}|\leq\frac{2\sqrt{1-\zeta}}{\sqrt{\left|\left(2(\gamma+2)\Upsilon^{k}_{3}C(\delta,3)+(\gamma-1)
    (\gamma+2)\left[\Upsilon^{k}_{2}C(\delta,2)\right]^{2}\right)\right|}}
\end{equation*}
and
\begin{equation*}\label{70}
    |a_{3}|\leq\frac{2(1-\zeta)}{(\gamma+2)\Upsilon^{k}_{3}C(\delta,3)}+\left[\frac{2(1-\zeta)}{(\gamma+1)\Upsilon^{k}_{2}C(\delta,2)}\right]^{2}.
\end{equation*}
\end{thm}

Putting $\gamma = 0$ in Theorem \ref{thm3.2}, we have the following corollary.

\begin{cor}
\textit{
Let the function $f(z)$ given by (\ref{1}) be in the class $B^{k,\alpha,\beta,\delta,\lambda}_{\Sigma}(\zeta)$. Then}
\begin{equation*}\label{71}
     |a_{2}|\leq\frac{\sqrt{2(1-\zeta)}}{\sqrt{\left|\left(2\Upsilon^{k}_{3}C(\delta,3)-\left[\Upsilon^{k}_{2}C(\delta,2)\right]^{2}\right)\right|}}
\end{equation*}
and
\begin{equation*}\label{72}
    |a_{3}|\leq\frac{(1-\zeta)}{\Upsilon^{k}_{3}C(\delta,3)}+\left[\frac{2(1-\zeta)}{\Upsilon^{k}_{2}C(\delta,2)}\right]^{2}.
\end{equation*}
\end{cor}
Putting $\gamma = 1$ in Theorem \ref{thm3.2}, we have the following corollary.

\begin{cor}
\textit{
Let the function $f(z)$ given by (\ref{1}) be in the class $B^{k,\alpha,\beta,\delta,\lambda}_{\Sigma}(1,\zeta)$. Then}
\begin{equation*}\label{73}
     |a_{2}|\leq\sqrt{\frac{2(1-\zeta)}{3\Upsilon^{k}_{3}C(\delta,3)}}
\end{equation*}
and
\begin{equation*}\label{74}
    |a_{3}|\leq\frac{2(1-\zeta)}{3\Upsilon^{k}_{3}C(\delta,3)}+\left[\frac{(1-\zeta)}{\Upsilon^{k}_{2}C(\delta,2)}\right]^{2}.
\end{equation*}
\end{cor}
Putting $k = 0$, from Corollaries 3.11 and 3.12, we get the following corollaries.

\begin{cor}

\textit{
Let the function $f(z)$ given by (\ref{1}) be in the class $B^{\alpha,\beta,\delta,\lambda}_{\Sigma}(\zeta)$. Then}

\begin{equation*}\label{75}
     |a_{2}|\leq\frac{2(1-\zeta)}{\sqrt{\left|\left(2C(\delta,3)-\left[C(\delta,2)\right]^{2}\right)\right|}}
\end{equation*}
and
\begin{equation*}\label{76}
    |a_{3}|\leq\frac{1-\zeta}{C(\delta,3)}+\left[\frac{(2(1-\zeta))}{C(\delta,2)}\right]^{2}.
\end{equation*}
\end{cor}
\begin{cor}
\textit{
Let the function $f(z)$ given by (\ref{1}) be in the class $B^{\alpha,\beta,\delta,\lambda}_{\Sigma}(1,\zeta)$. Then}

\begin{equation*}\label{77}
     |a_{2}|\leq\sqrt{\frac{2(1-\zeta)}{3C(\delta,3)}}
\end{equation*}
and
\begin{equation*}\label{78}
    |a_{3}|\leq\frac{2(1-\zeta)}{3C(\delta,3)}+\left[\frac{(1-\zeta)}{C(\delta,2)}\right]^{2}.
\end{equation*}
\end{cor}

Putting $\delta = 0$, from Corollaries 3.11 and 3.12, we get the following corollaries.

\begin{cor}
\textit{
Let the function $f(z)$ given by (\ref{1}) be in the class $B^{k,\alpha,\beta,\lambda}_{\Sigma}(\zeta)$. Then}
\begin{equation*}\label{79}
     |a_{2}|\leq\frac{\sqrt{2(1-\zeta)}}{\sqrt{\left|\left(2\Upsilon^{k}_{3}-[\Upsilon^{k}_{2}]^{2}\right)\right|}}
\end{equation*}
and
\begin{equation*}\label{80}
    |a_{3}|\leq\frac{(1-\zeta)}{\Upsilon^{k}_{3}}+\left[\frac{2(1-\zeta)}{\Upsilon^{k}_{2}}\right]^{2}.
\end{equation*}
\end{cor}
\begin{cor}
\textit{
Let the function $f(z)$ given by (\ref{1}) be in the class $B^{k,\alpha,\beta,\lambda}_{\Sigma}(1,\zeta)$. Then}
\begin{equation*}\label{81}
     |a_{2}|\leq\sqrt{\frac{2(1-\zeta)}{3\Upsilon^{k}_{3}}}
\end{equation*}
and
\begin{equation*}\label{82}
    |a_{3}|\leq\frac{2(1-\zeta)}{3\Upsilon^{k}_{3}}+\left[\frac{(1-\zeta)}{\Upsilon^{k}_{2}}\right]^{2}.
\end{equation*}
\end{cor}

Putting $k=\delta=0$, from Corollaries 3.11 and 3.12, we get the following corollaries.
\begin{cor}
\textit{
Let the function $f(z)$ given by (\ref{1}) be in the class $B^{\alpha,\beta,\lambda}_{\Sigma}(\zeta)$. Then}
\begin{equation*}\label{83}
     |a_{2}|\leq\sqrt{2(1-\zeta)}
\end{equation*}
and
\begin{equation*}\label{84}
    |a_{3}|\leq(1-\zeta)+4\left[(1-\zeta)\right]^{2}.
\end{equation*}
\end{cor}

\begin{cor}
\textit{
Let the function $f(z)$ given by (\ref{1}) be in the class $B^{\alpha,\beta,\lambda}_{\Sigma}(1,\zeta)$. Then}
\begin{equation*}\label{85}
     |a_{2}|\leq\sqrt{\frac{2(1-\zeta)}{3}}
\end{equation*}
and
\begin{equation*}\label{86}
    |a_{3}|\leq\frac{2(1-\zeta)}{3}+\left[1-\zeta\right]^{2}.
\end{equation*}
\end{cor}

\textbf{Concluding Remarks}: By specializing the parameters of operator, various other interesting corollaries and consequences of our main results (which are asserted by Theorem \ref{thm2.1} above)  can be derived easily. The details involved shall be left as an exercise for the interested readers.

\medskip\noindent


\bigskip
\noindent
Adnan Ghazy Alamoush\\
Faculty of Science, Taibah University, Saudi Aarabia.
Email:adnan\_omoush@yahoo.com\\
\end{document}